\documentclass[aoas,dvips]{arximspdf}
\usepackage{dcolumn}
\usepackage{graphicx}

\doi{10.1214/09-AOAS312A}
\referstodoi{10.1214/09-AOAS312}
\volume{3}
\issue{4}
\pubyear{2009}
\firstpage{1266}
\lastpage{1269}

\makeatletter
\newcolumntype{d}[1]{D{.}{.}{#1}}
\def\R{\mathbb{R}}
\def\E{\mathbb{E}}
\def\P{\mathbb{P}}

\makeatother

\begin{document}
\begin{frontmatter}

\title{Discussion of: Brownian distance covariance}
\pdftitle{Discussion on Brownian distance covariance by G. J. Szekely and M. L. Rizzo}
\runtitle{Discussion}

\begin{aug}
\author[a]{\fnms{Peter J.} \snm{Bickel}\corref{}\ead[label=e1]{bickel@stat.berkeley.edu}\thanksref{r1}}
\and
\author[a]{\fnms{Ying} \snm{Xu}\ead[label=e2]{yingxu@stat.berkeley.edu}\thanksref{r1}}
\runauthor{P. J. Bickel and Y. Xu}
\affiliation{University of
California,  Berkeley}
\address[a]{Department of Statistics\\
367 Evans Hall\\
Berkeley, California 94710--3860\\
USA\\
\printead{e1}\\
\phantom{\textsc{E-mail}: }\printead*{e2}} %adresu isvedimo komanda gale!
\end{aug}
\thankstext{r1}{Supported in part by NSF Grant DMS-09-06808.}
% HISTORY:

% ABSTRACT

% KEYWORDS

\end{frontmatter}

Szekely and Rizzo present a new interesting measure of correlation.
The idea of using $\int|\phi_n(u,v)-\phi_n^{(1)}(u)\phi
_n^{(2)}(v)|^2\,d\mu(u,v)$, where $\phi_n$, $\phi_n^{(1)}$, $\phi
_n^{(2)}$ are the empirical characteristic functions of a sample
$(X_i,Y_i)$, $i=1,\ldots, n$, of independent copies of $X$ and $Y$ is
not so novel. A. Feuerverger considered such measures in a series of
papers \cite{4}. Aiyou Chen and I have actually analyzed such a measure for
estimation in \cite{3} in connection with ICA.

However, the choice of $\mu(\cdot, \cdot)$ which makes the measure
scale free, the extension to $X\in\R^p,Y \in\R^q$ and its
identification with the Brownian distance covariance is new, surprising
and interesting.

There are three other measures available, for general
$p$, $q$:

\begin{enumerate}
\item The canonical correlation $\rho$ between $X$ and $Y$.
\item The rank correlation $r$ (for $p=q=1$) and its canonical
correlation generalization.
\item The Renyi correlation $R$.
\end{enumerate}

All vanish along with the Brownian distance (BD) correlation in the
case of independence and all are scale free. The Brownian distance and
Renyi covariance are the only ones which vanish iff $X$ and $Y$ are
independent.

However, the three classical measures also give a characterization of
total dependence. If $|\rho|=1$, $X$ and $Y$ must be linearly related;
if $|r|=1$, $Y$ must be a monotone function of $X$ and if $R=1$, then
either there exist nontrivial functions~$f$ and $g$ such that $\P
(f(X)=g(Y))=1$ or at least there is a sequence of such nontrivial
functions $f_n$, $g_n$ of variance $1$ such that $\E
(f_n(X)-g_n(Y))^2\rightarrow0$.

In this respect, by Theorem 4 of Szekely and Rizzo, for the common
$p=q=1$ case, BD correlation does not differ from Pearson correlation.

Although we found the examples varied and interesting and the
computation of~$p$ values for the BD covariance effective, we are not
convinced that the comparison with the rank and Pearson correlations is
quite fair, and think a comparison to~$R$ is illuminating.

Intuitively, the closer the form of observed dependence is to that
exhibited for the extremal value of the statistic, the more power one
should expect. Example 1 has $Y$ as a distinctly nonmonotone function
of $X$ plus noise, a situation where we would expect the rank
correlation to be weak and, similarly, the other examples correspond to
nonlinear relationships between $X$ and $Y$ in which we would expect
the Pearson correlation to perform badly.
In general, for goodness of fit, it is important to have statistics
with power in directions which are plausible departures; see Bickel, Ritov  and  Stoker \cite{1}.

Ying Xu is studying, in the context of high dimensional data, a version
of empirical Renyi correlation different from that of Breiman and
Friedman \cite{2}.

Let $f_1, f_2, \ldots$ be an orthonormal basis of $L_2(P_X)$ and $g_1,
g_2, \ldots$ an orthonormal basis of $L_2(P_Y)$, where $L_2(P_X)$ is
the Hilbert space of function $f$ such that $\E f^2(X)<\infty$ and
similarly for $L_2(P_Y)$.

Let the $(K,L)$ approximate Renyi correlation be defined as
\[
\max\Biggl\{ \mbox{corr}\Biggl(\sum_{k=1}^K\alpha_kf_k(X),\sum_{l=1}^L\beta_l
g_l(Y) \Biggr) \Biggr\},
\]
where corr is Pearson correlation.

This is seen to be the canonical correlation of $\underline{f}(X)$ and
$\underline{g}(Y)$, where $\underline{f}\equiv(f_1, \ldots,
f_K)^T$, $\underline{g}\equiv(g_1, \ldots, g_L)^T$, and is easily
calculated as a generalized eigenvalue problem. The empirical $(K,L)$
correlation is just the solution of the corresponding empirical problem
where the variance covariance matrices $\operatorname{Var} \underline{f}(X)\equiv\E
[\underline{f_c}(X)\underline{f_c}^T(X)]$ where $ \underline
{f_c}(X)\equiv \underline{f}(X)-\E\underline{f}(X) $,
$\operatorname{Var} \underline{g}(Y)$ and $\operatorname{Cov}( \underline{f}(X),\underline
{g}(Y))$ are replaced by their empirical counterparts. For
$K,L\rightarrow\infty$, the $(K,L)$ correlation tends to the Renyi
correlation,
\[
R\equiv\max\{ \mbox{corr}( f(X), g(Y) ) \dvtx f\in L_2(P_X), g \in
L_2(P_Y) \}.
\]
For the empirical $(K,L)$ correlation, $K$ and $L$ have to be chosen in
a data determined way, although evidently each $K$, $L$ pair provides a
test statistic. An even more important choice is that of the $f_k$ and
$g_l$ (which need not be orthonormal but need only have a linear span
dense in their corresponding Hilbert spaces).

We compare the performance of these test statistics in the first of the
Szekely--Rizzo examples in the next section.

\section{Comparison on data example}

Here we will investigate the performance of the standard ACE estimate of
the Renyi correlation and a version of $(K,L)$ correlation in the first
of the Szekely--Rizzo examples.

Breiman and Friedman \cite{2} provided an algorithm, known as alternating
conditional expectations (ACE), for estimating the transformations
$f_0$, $g_0$ and $R$ itself.

\begin{figure}

\includegraphics{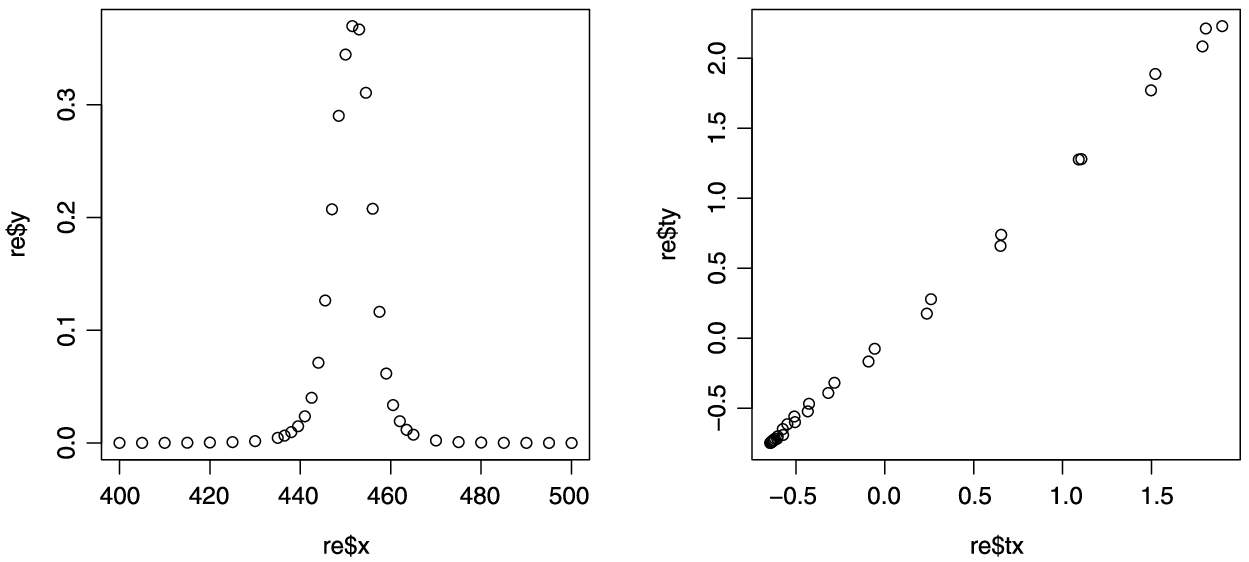}

   \caption{}
  \label{fig1}
\end{figure}

The estimated Renyi correlation is very close to $1$ ($0.9992669$) in
this case, as expected
since $Y$ is a function of $X$ plus some noise. Figure~\ref{fig1} shows
the original relationship between $X$ and $Y$ on the left and the
relationship between the estimated transformations $\hat{f}$ and $\hat
{g}$ on the right.
Having computed $\hat{R}$, the estimate of~$R$, we compute its
significance under the null hypothesis of independence using the
permutation distribution just as Szekely and Rizzo did. The $p$-value is
${\leq}0.001$, which is extremely small as it should be.

\begin{table}
\caption{}
\label{tab1}
\begin{tabular*}{\textwidth}{@{\extracolsep{\fill}}ld{1.7}d{1.7}d{2.6}@{}}
\hline
&\multicolumn{1}{c}{$\bolds{K=2}$, $\bolds{L=2}$}  & \multicolumn{1}{c}{$\bolds{K=3}$, $\bolds{L=4}$}
& \multicolumn{1}{c@{}}{$\bolds{K=5}$, $\bolds{L=5}$}\\
\hline
Estimated ($K$,$L$) correlation & 0.8160803 &0.9170764 & 0.977163\\
$p$-value &0.002 & 0.002 & {\leq}0.001\\
\hline
\end{tabular*}
\end{table}

Next, we compute the empirical $(K,L)$ correlation. Given that the
proposed nonlinear model is
\[
y=\frac{\beta_1}{\beta_2}\exp\biggl\{\frac{-(x-\beta_3)^2}{2\beta
_2^2}\biggr\}+\varepsilon,
\]
we chose, as an orthonormal basis with respect to the Lebesgue measure,
one defined by the Hermite polynomials defined as $H_n(x)=(-1)^n
e^{x^2/2}\frac{d^n}{dx^n}e^{-x^2/2} $, for both $X$ and $Y$. We take
$f_k(\cdot)=g_k(\cdot)=e^{-x^2/4}H_k(\cdot)$.

Table~\ref{tab1} gives the computation results of different
combinations of $K$ and $L$. As before, the $p$-value is computed by a
permutation test, based on $999$ replicates.

The value, not surprisingly, is close to $\hat{R}$, for $K=L=5$.

\printaddresses

\end{document}